\documentclass{amsart}

\usepackage{latexsym,enumerate}
\usepackage{amsmath,amsthm,amsopn,amstext,amscd,amsfonts,amssymb}
\usepackage[ansinew]{inputenc}
\usepackage{verbatim}
\usepackage{graphicx}
\usepackage{pstricks}
\usepackage{wasysym}
\usepackage[all]{xy}
\usepackage{epsfig} 
\usepackage{graphicx,psfrag} 
\usepackage{subfigure}
\usepackage{pstricks,pst-node}

\usepackage{multicol}
\usepackage{color}
\usepackage{colortbl}

\newtheorem{theorem}{Theorem}[section]
\newtheorem{lemma}[theorem]{Lemma}
\newtheorem{proposition}[theorem]{Proposition}
\newtheorem{corollary}[theorem]{Corollary}

\newtheorem{example}[theorem]{Example}
\newtheorem{remark}[theorem]{Remark}

\newcommand{\spb}[1]{\smallskip}
\newcommand{\mpb}[1]{\medskip}
\newcommand{\bpb}[1]{\bigskip}

\renewcommand{\d}{\delta}
\newcommand{\D}{\Delta}

\newcommand{\G}{\Gamma}

\begin{document}
\DeclareGraphicsExtensions{.jpg,.pdf,.mps,.png}

\title{New lower bounds for the Geometric-Arithmetic index}

\author[Alvaro Mart\'{\i}nez-P\'erez]{Alvaro Mart\'{\i}nez-P\'erez$^{(1)}$}
\address{ Facultad CC. Sociales de Talavera,
Avda. Real Fábrica de Seda, s/n. 45600 Talavera de la Reina, Toledo, Spain}
\email{alvaro.martinezperez@uclm.es}
\thanks{$^{(1)}$ Supported in part by a grant
from Ministerio de Econom{\'\i}a y Competitividad (MTM 2015-63612P), Spain.}

\author[Jos\'e M. Rodr{\'\i}guez]{Jos\'e M. Rodr{\'\i}guez$^{(2)}$}
\address{Departamento de Matem\'aticas, Universidad Carlos III de Madrid,
Avenida de la Universidad 30, 28911 Legan\'es, Madrid, Spain}
\email{jomaro@math.uc3m.es}
\thanks{$^{(2)}$ Supported in part by two grants
from Ministerio de Econom{\'\i}a y Competitividad, Agencia Estatal de
Investigación (AEI) and Fondo Europeo de Desarrollo Regional (FEDER) (MTM 2016-78227-C2-1-P and MTM 2015-69323-REDT), Spain, and a grant from CONACYT (FOMIX-CONACyT-UAGro 249818), M\'exico.}

\date{\today}

\begin{abstract}
The concept of geometric-arithmetic index was introduced in the
chemical graph theory recently, but it has shown to be useful.
The aim of this paper is to obtain new inequalities involving the geometric-arithmetic index $GA_1$
and characterize graphs extremal with respect to them.
Our main results provide lower bounds $GA_1(G)$ involving just the minimum and the maximum degree of the graph $G$.
\end{abstract}

\maketitle{}

{\it Keywords:  Geometric-arithmetic index, Graph invariant, Vertex-degree-based graph invariant, Topological index.}

{\it 2010 AMS Subject Classification numbers: 05C07, 92E10.} 

\section{Introduction}
A single number, representing a chemical structure in graph-theoretical terms via the
molecular graph, is called a topological descriptor and if it in addition correlates with
a molecular property it is called topological index, which is used to understand physicochemical
properties of chemical compounds.
Topological indices are interesting since they capture some of the properties of a molecule in a single number.
Hundreds of topological indices have been introduced and studied, starting with the
seminal work by Wiener in which he used the sum of all shortest-path distances of
a (molecular) graph for modeling physical properties of alkanes (see \cite{Wi}).

Topological indices based on end-vertex degrees of edges have been
used over 40 years. Among them, several indices are recognized to be useful tools in
chemical researches. Probably, the best know such descriptor is the Randi\'c connectivity
index ($R$) \cite{R}. There are more than thousand papers and a couple of books dealing with
this molecular descriptor (see, e.g., \cite{GF}, \cite{LG}, \cite{LS}, \cite{RS}, \cite{RS0} and the references therein).
During many years, scientists were trying to improve the predictive power of the
Randi\'c index. This led to the introduction of a large number of new topological
descriptors resembling the original Randi\'c index.
The first geometric-arithmetic index $GA_1$, defined in \cite{VF} as
$$
GA_1 = GA_1(G) = \sum_{uv\in E(G)}\frac{\sqrt{d_u d_v}}{\frac12 (d_u + d_v)}
$$
where $uv$ denotes the edge of the graph $G$ connecting the vertices $u$ and $v$, and
$d_u$ is the degree of the vertex $u$,
is one of the successors of the Randi\'c index.
Although $GA_1$ was introduced in $2009$, there are many papers dealing with this index
(see, e.g., \cite{Das10b}, \cite{DGF}, \cite{DGF2}, \cite{MH}, \cite{RS3}, \cite{S}, \cite{VF} and the references therein).
There are other geometric-arithmetic indices, like $Z_{p,q}$ ($Z_{0,1} = GA_1$), but the results in \cite[p.598]{DGF}
show that the $GA_1$ index gathers the
same information on observed molecule as other $Z_{p,q}$ indices.

The reason for introducing a new index is to gain prediction of target property (properties)
of molecules somewhat better than obtained by already presented indices. Therefore,
a test study of predictive power of a new index must be done. As a standard for
testing new topological descriptors, the properties of octanes are commonly used.
We can find 16 physico-chemical properties of octanes at www.moleculardescriptors.eu.

The $GA_1$ index gives better correlation coefficients than $R$ for these properties, but the differences between
them are not significant. However, the predicting ability of the $GA_1$ index compared with
Randi\'c index is reasonably better (see \cite[Table 1]{DGF}).
Although only about 1000 benzenoid hydrocarbons are known, the number of
possible benzenoid hydrocarbons is huge. For instance, the number of
possible benzenoid hydrocarbons with 35 benzene rings is $5.85\cdot 10^{21}$ \cite{NGJ}.
Therefore, modeling their physico-chemical properties is important in order
to predict properties of currently unknown species.
The graphic in \cite[Fig.7]{DGF} (from \cite[Table 2]{DGF}, \cite{TRC}) shows
that there exists a good linear correlation between $GA_1$ and the heat of formation of benzenoid hydrocarbons
(the correlation coefficient is equal to $0.972$).

Furthermore, the improvement in
prediction with $GA_1$ index comparing to Randi\'c index in the case of standard
enthalpy of vaporization is more than 9$\%$. That is why one can think that $GA_1$ index
should be considered in the QSPR/QSAR researches.

Throughout this work, $G=(V (G),E (G))$ denotes a (nonoriented) finite simple (without multiple edges and loops) nontrivial ($E(G) \neq \emptyset$) graph.
The aim of this paper is to obtain new inequalities involving the geometric-arithmetic index $GA_1$
and characterize graphs extremal with respect to them.
Our main results provide lower bounds $GA_1(G)$ involving just the minimum and the maximum degree of the graph $G$ (see Theorems \ref{Th: lower bound} and \ref{l:230}, and Corollaries \ref{C:d-hd 2} and \ref{C:d-1d}).

\section{$GA_1$ and minimum and maximum degree}

If $G$ is a graph with $m$ edges, minimum degree $\delta$ and maximum degree $\Delta$, then in \cite{Das10b} (see also \cite{DGF}) we find the bounds:
\begin{equation}\label{eq: bound}
\frac{2m\sqrt{\delta \Delta}}{\delta+\Delta} \leq GA_1(G)\leq m.
\end{equation}

\begin{remark} \label{remark bipartite}
$GA_1(G)=\frac{2m\sqrt{\delta \Delta}}{\delta+\Delta}$ if and only if the graph is
either regular or bipartite with the two sets being respectively the set of vertices with degree $\delta$ and degree $\Delta$.
\end{remark}

Let us recall Lemma 2.2 and Corollary 2.3 in \cite{RS2}.

\begin{lemma}\label{lema 1} Let $f$ be the function $f(t) = \frac{2t}{1+t^2}$ on the interval $[0,\infty)$. Then $f$ strictly
increases in $[0, 1]$, strictly decreases in $[1,\infty)$, $f(t) = 1$ if and only if $t = 1$ and $f(t) = f(t_0)$
if and only if either $t = t_0$ or $t = t_0^{-1}$.
\end{lemma}

\begin{corollary}\label{cor} Let $g$ be the function $g(x, y) = \frac{2\sqrt{xy}}{x+y}$ with $0 < a\leq  x, y \leq b$. Then
$\frac{2\sqrt{ab}}{a+b}\leq g(x,y)\leq 1.$
The equality in the lower bound is attained if and only if either $x = a$ and $y = b$, or $x = b$
and $y = a$, and the equality in the upper bound is attained if and only if $x = y$.
\end{corollary}

Given integers $0 < \delta \le \Delta$, let us define $\mathcal{G}_{\d,\D}$ as the set of graphs $G$ with minimum degree $\delta$, maximum degree $\Delta$ and such that:

$(1)$ $G$ is isomorphic to the complete graph with $\D+1$ vertices $K_{\D+1}$, if $\d = \D$,

$(2)$ $|V(G)|=\D+1$, there are $\D$ vertices with degree $\d$, if $\d < \D$ and $\Delta (\delta+1)$ is even,

$(3)$ $|V(G)|=\D+1$, there are $\D-1$ vertices with degree $\d$ and a vertex with degree $\d+1$, if $\d < \D-1$ and $\Delta (\delta+1)$ is odd,

$(4)$ $|V(G)|=\D+1$, there are $\D-1$ vertices with degree $\d$ and two vertices with degree $\D$, if $\d = \D-1$ and $\Delta$ is odd (and thus $\Delta (\delta+1)$ is odd).

\begin{remark} \label{remark sharp0}
Every graph $G \in \mathcal{G}_{\d,\D}$ has maximum degree $\Delta$ and $|V(G)|=\D+1$.
Hence, every graph $G \in \mathcal{G}_{\d,\D}$ is connected.
\end{remark}

\begin{proposition}\label{prop edges}
For any integers $0 < \delta \le \Delta$, we have $\mathcal{G}_{\d,\D} \neq \emptyset$.
Let $G$ be a graph with minimum degree $\delta$ and maximum degree $\Delta$. Then
\[|E(G)|\geq \frac{\Delta(\delta+1)}{2}  \quad \mbox{if} \ \Delta (\delta+1) \ \mbox{is even},\quad
|E(G)|\geq \frac{\Delta(\delta+1)+1}{2}  \quad \mbox{if} \ \Delta  (\delta+1) \ \mbox{is odd},\]
with equality if and only if $G \in \mathcal{G}_{\d,\D}$.
\end{proposition}

\begin{proof} There is at least one vertex $v_0 \in V(G)$ with degree $\Delta$ and $\Delta$ vertices, $v_1,\dots,v_\Delta$, adjacent to it.
Since $d_{v_i}\geq \delta$, $|E(G)|\geq \frac{\Delta + \Delta\delta}{2}=\frac{\Delta(\delta+1)}{2}$.
If $\Delta (\delta+1)$ is odd, then $\frac{\Delta(\delta+1)}{2}$ is not an integer and the lower bound is at least $\frac{\Delta(\delta+1)+1}{2}$.

If the equality is attained, then $V(G)=\{v_0,v_1,\dots,v_\Delta\}$ (thus $|V(G)|=\D+1$). As in Remark \ref{remark sharp0}, we can conclude that $G$ is connected.
If $\Delta (\delta+1)$ is even, then $d_{v_i} = \delta$ for $i=1,\dots,\D,$ and $G \in \mathcal{G}_{\d,\D}$.
If $G \in \mathcal{G}_{\d,\D}$, it is clear that the equality holds.
If $\Delta (\delta+1)$ is odd, then a similar argument gives that the equality is attained if and only if
$G \in \mathcal{G}_{\d,\D}$.

Finally, let us prove $\mathcal{G}_{\d,\D} \neq \emptyset$.
This is clear in the case $(1)$; so, let us assume $\d < \D$.
Consider a graph $H$ with $\Delta+1$ vertices, $v_0,v_1,\dots,v_\Delta$.
Assume that $d_{v_0}=\Delta$. Then, there is an edge joining $v_0$ with $v_i$ for every $i>0$.

First, suppose $\delta$ is odd.
Thus, $\Delta (\delta+1)$ is even.
We have already one edge in each $v_i$.
We are going to add edges so that $d_{v_i}=\delta$ for every $i>0$.
Let us define for every $1 \le i,j\le \D$, $||i-j||=\min\{|i-j|,\Delta-|i-j|\}$ (this is, the distance between the vertices $v_i$ and $v_j$ in the cycle $v_1,v_2,\dots,v_{\Delta},v_1$).
Consider an edge $v_iv_j$ for every pair of vertices with $||i-j||\leq \frac{\delta-1}{2}$.
This is possible since $\delta-1$ is even and $\delta-1<\Delta-1$. Then, every vertex $v_i$ with $i>0$ satisfies that $\d_{v_i}=\delta$ and $H \in \mathcal{G}_{\d,\D}$.

Now, suppose $\delta$ and $\Delta$ are even.
Thus, $\Delta (\delta+1)$ is even.
Consider an edge $v_iv_j$ for every pair of vertices with $||i-j||\leq \frac{\delta-2}{2} $ and a edge $v_iv_j$ for every $||i-j||=\frac{\Delta}{2}$. Notice that this is well defined since $\Delta$ is even and it is a  new edge since $\frac{\Delta}{2}>\frac{\delta-2}{2}$. Then, every vertex $v_i$ with $i>0$ satisfies that $d_{v_i}=\delta$ and $H \in \mathcal{G}_{\d,\D}$.

Finally, if $\delta$ is even and $\Delta$ is odd, then $\Delta (\delta+1)$ is odd.
Consider an edge $v_iv_j$ for every pair of vertices with $||i-j||\leq \frac{\delta-2}{2}$.
Now every vertex $v_i$ with $i>0$ has degree $\delta-1$.
Let us define, for every $1\leq i<j\leq \Delta-1$, an edge $v_iv_j$ if $j-i=\frac{\Delta-1}{2}$. This edge is new since $\frac{\delta-2}{2}<\frac{\Delta-1}{2}$. Now, $d_{v_i}=\delta$ for every $0<i< \Delta$. It suffices to define an edge joining $v_{\Delta}$ to any non-adjacent vertex $v_{i_0}$, for example $i_0=\frac{\delta}{2}+1$, and therefore, $H \in \mathcal{G}_{\d,\D}$.
Notice that, in this case, $d_{v_0}=\Delta$, $d_{v_{i_0}}=\delta+1$ and $d_{v_i}=\delta$ for every $i\neq 0,i_0$.
\end{proof}

\begin{proposition}\label{Prop minimal}
For every integers $0< \delta \le \Delta$ and $G \in \mathcal{G}_{\d,\D}$, we have
$$
\begin{aligned}
GA_1(G)
& = \frac{2\Delta \sqrt{\delta \Delta}}{\delta+\Delta}+\frac{\Delta(\delta-1)}{2}  \quad \mbox{if } \ \Delta (\delta+1) \ \mbox{is even},
\\
GA_1(G) &
= \frac{2(\Delta-1) \sqrt{\delta \Delta}}{\delta+\Delta}+\frac{2 \sqrt{(\delta+1) \Delta}}{\delta+1+\Delta}+\frac{2\delta \sqrt{\delta (\delta+1)}}{2\delta+1}+\frac{(\Delta-2)(\delta-1)-1}{2} \\
& > \frac{2\Delta \sqrt{\delta \Delta}}{\delta+\Delta}+\frac{\Delta(\delta-1)}{2} \quad \mbox{if } \ \Delta (\delta+1) \ \mbox{is odd}.
\end{aligned}
$$
\end{proposition}

\begin{proof} The equalities follow from the definitions of $GA_1$ and $\mathcal{G}_{\d,\D}$. Let us see that
\[\frac{2(\Delta-1) \sqrt{\delta \Delta}}{\delta+\Delta}+\frac{2 \sqrt{(\delta+1) \Delta}}{\delta+1+\Delta}+\frac{2\delta \sqrt{\delta (\delta+1)}}{2\delta+1}+\frac{(\Delta-2)(\delta-1)-1}{2}  > \frac{2\Delta \sqrt{\delta \Delta}}{\delta+\Delta}+\frac{\Delta(\delta-1)}{2}\]
if $\Delta (\delta+1)$ is odd.
If $\delta=\Delta$, then $\Delta(\delta+1)$ is even.
Thus, we can assume that $\delta<\Delta$.
It suffices to check that
$$
\frac{2\sqrt{(\delta+1) \Delta}}{\delta+1+\Delta} >\frac{2 \sqrt{\delta \Delta}}{\delta+\Delta}
\qquad
\text{and}
\qquad
\frac{2\delta \sqrt{\delta (\delta+1)}}{2\delta+1}-\frac{1}{2}>\delta-1.
$$
The first claim follows from Lemma \ref{lema 1} and the fact that $\frac{2\sqrt{xy}}{x+y}= f(t)$ with $t = \sqrt{\frac{x}{y}}$, since $1\leq \frac{\Delta}{\delta+1}<\frac{\Delta}{\delta}$.
For the second claim it suffices to check that
$$\frac{2\delta \sqrt{\delta (\delta+1)}}{2\delta+1}>\frac{2\delta-1}{2} \quad \Leftrightarrow \quad \sqrt{\delta(\delta+1)}>\frac{4\delta^2-1}{4\delta} \quad \Leftrightarrow$$
$$\delta(\delta+1)>\frac{16\delta^4-8\delta^2+1}{16\delta^2} \quad \Leftrightarrow \quad 16\delta^4+16\delta^3>16\delta^4-8\delta^2+1 \quad \Leftrightarrow \quad 16\delta^3+8\delta^2>1,$$
finishing the proof.
\end{proof}

\begin{theorem} \label{Th: lower bound}
Let $G$ be a graph with minimum degree $\delta>0$ and maximum degree $\Delta$. Then
$$
GA_1(G)\geq \Delta(\delta+1)\frac{\sqrt{\delta\Delta}}{\delta+\Delta}
$$
with equality if and only if $\d=1$ and $G$ is a star graph or $\d=\D$ and $G$ is a complete graph.

\smallskip

Furthermore, if $\Delta(\delta+1)$ is odd, then
$$
GA_1(G)\geq \big(\Delta(\delta+1)+1\big)\frac{\sqrt{\delta\Delta}}{\delta+\Delta}.
$$
\end{theorem}

\begin{proof}
By Proposition \ref{prop edges}, $2m\geq \Delta(\delta+1)$.
This inequality and the lower bound in \eqref{eq: bound} give the first inequality.
If we use in this argument the second part of Proposition \ref{prop edges}, then we obtain the second inequality.

If $\d=1$ and $G$ is a star graph or $\d=\D$ and $G$ is a complete graph, then one can check that the equality is attained in the first inequality.

If the equality holds in the first inequality for a graph $G$, then Remark \ref{remark bipartite} gives that $G$ is
either regular or bipartite with the two sets being respectively the set of vertices with degree $\delta$ and degree $\Delta$,
and Proposition \ref{prop edges} gives that $\Delta(\delta+1)$ is even and $G \in \mathcal{G}_{\d,\D}$.
If $\d=\D$, then $G$ is a complete graph.
If $\d<\D$, then $G$ is a bipartite graph with the two sets being the set of vertices with degree $\delta$ and degree $\Delta$;
thus, given a vertex $v$ with degree $\d$, there are $\d-1$ edges connecting $v$ and $\d-1$ vertices with degree $\d$;
if $\d>1$, then this is not possible since $G$ is a bipartite graph. Hence, $\d = 1$ and $G$ is a star graph.
\end{proof}

We say that a graph $G$ with minimum degree $\delta$ and maximum degree $\Delta$ is \emph{minimal} if $GA_1(G) \le GA_1(\G)$
for every graph $\G$ with minimum degree $\delta$ and maximum degree $\Delta$.

\begin{proposition}\label{prop vertex-edges}
For any integers $0 < \delta \le \Delta$, let $G$ be a graph with minimum degree $\delta$ and maximum degree $\Delta$ which is minimal for those $\delta$ and $\Delta$. Then
\[ 
\begin{aligned}
\frac{\Delta(\delta+1)}{2} 
& \leq |E(G)|\leq \Delta \delta  \quad \mbox{if} \ \Delta (\delta+1) \ \mbox{is even},
\\
\frac{\Delta(\delta+1)+1}{2}
& \leq |E(G)|\leq \Delta \delta \quad \mbox{if} \ \Delta (\delta+1) \ \mbox{is odd},
\\
\Delta+1
& \leq |V(G)|\leq \frac{\Delta(2\delta-1)}{\delta}+1.
\end{aligned}
\]
For $|E(G)|$, the lower bound  is attained if and only if $G\in \mathcal{G}_{\delta,\Delta}$ and the upper bound is attained if and only if
$G=K_{\delta,\Delta}$.
\end{proposition}

\begin{proof} By Proposition \ref{prop edges}, if $\Delta (\delta+1)$ is even, $\frac{\Delta(\delta+1)}{2} \leq |E(G)|$ and, if $\Delta (\delta+1)$ is odd, $\frac{\Delta(\delta+1)+1}{2} \leq |E(G)|$. By Corollary \ref{cor}, for every edge $uv$,
$\frac{2\sqrt{\delta \Delta}}{\delta+\Delta}\leq \frac{2\sqrt{d_ud_v}}{d_u+d_v}.$ Therefore, by (\ref{eq: bound}), if $|E(G)|>\Delta \delta$, then $GA_1(G)>\Delta \delta \frac{2\sqrt{\delta \Delta}}{\delta+\Delta}=GA_1(K_{\delta,\Delta})$ leading to contradiction.

It is immediate to see that $\Delta+1 \leq |V(G)|$ since there is a vertex with degree $\Delta$ and $\Delta$ vertices adjacent to it. Now suppose $|V(G)|>\frac{\Delta(2\delta-1)}{\delta}+1$. By hypothesis, there is a vertex with degree $\Delta$ and more than $\frac{\Delta(2\delta-1)}{\delta}$ vertices with degree at least $\delta$. Thus, $|E(G)|>\delta\Delta$, leading to contradiction.

By Proposition \ref{prop edges}, the lower bound for $|E(G)|$ is attained if and only if $G\in \mathcal{G}_{\delta,\Delta}$. By Corollary \ref{cor}, the upper bound for $|E(G)|$ is attained if and only if $G=K_{\delta,\Delta}$.
\end{proof}

Let us denote by $K_{\delta,\Delta}$ the complete bipartite graph with a partition $K_1$, $K_2$ with $\delta$ and $\Delta$ vertices respectively.
Notice that the vertices in $K_1$ have degree $\Delta$ and the vertices in $K_2$ have degree $\delta$.
It was proved in \cite{RS2} that
$GA_1(K_{\delta,\Delta})=\frac{2\delta \Delta \sqrt{\delta \Delta}}{\delta+\Delta}.$

Let $H_{\delta,\Delta}$ be any graph in $\mathcal{G}_{\delta,\Delta}$.
Note that if $\d=1$, then $H_{1,\Delta}=K_{1,\Delta}$.

\begin{proposition}\label{Prop: comp}
For any integers $1<\delta \leq \Delta$, we have

$(1)$ if $\frac{\Delta}{\delta}>\big( 2+\sqrt{3}\, \big)^2$, then $GA_1(H_{\delta,\Delta})>GA_1(K_{\delta,\Delta})$,

$(2)$ if $\frac{\Delta}{\delta}<\big( 2+\sqrt{3}\, \big)^2$ and $\D(\d+1)$ is even, then $GA_1(H_{\delta,\Delta})<GA_1(K_{\delta,\Delta})$.
\end{proposition}

\begin{proof} By Lemma \ref{lema 1}, $f(t)$ is decreasing in $[1,\infty)$ and $f\Big(\sqrt{\frac{\Delta}{\delta}}\;\Big)=\frac{2\sqrt{\delta \Delta}}{\delta+\Delta}$.

If $\frac{\Delta}{\delta}> \big( 2+\sqrt{3}\, \big)^2$, then
$$
\frac{2\sqrt{\delta \Delta}}{\delta+\Delta}
=f\left(\sqrt{\frac{\Delta}{\delta}}\;\right)
<f\big( 2+\sqrt{3}\, \big)
=\frac{2\big( 2+\sqrt{3}\, \big)}{1+\big( 2+\sqrt{3}\, \big)^2}
=\frac{1}{2}.
$$
Therefore, Proposition \ref{Prop minimal} and $\d>1$ give
\[
GA_1(H_{\delta,\Delta})
\ge \Delta\frac{2\sqrt{\delta \Delta}}{\delta+\Delta} + \frac{\Delta (\delta-1)}{2}
> \Delta \frac{2\sqrt{\delta \Delta}}{\delta+\Delta}
+ \Delta (\delta-1) \frac{2\sqrt{\delta \Delta}}{\delta+\Delta}
= GA_1(K_{\delta,\Delta}).
\]

If $\frac{\Delta}{\delta}< \big( 2+\sqrt{3}\, \big)^2$ and $\D(\d+1)$ is even, then
$$
\begin{aligned}
\frac{2\sqrt{\delta \Delta}}{\delta+\Delta}
& =f\left(\sqrt{\frac{\Delta}{\delta}}\;\right)
> f\big( 2+\sqrt{3}\, \big)
=\frac{2\big( 2+\sqrt{3}\, \big)}{1+\big( 2+\sqrt{3}\, \big)^2}
=\frac{1}{2} ,
\\
GA_1(H_{\delta,\Delta})
& =\Delta\frac{2\sqrt{\delta \Delta}}{\delta+\Delta} + \frac{\Delta (\delta-1)}{2}
<\Delta \delta \frac{2\sqrt{\delta \Delta}}{\delta+\Delta}
=GA_1(K_{\delta,\Delta}).
\end{aligned}
$$
\end{proof}


It may be wondered if
\begin{equation} \label{eq:conj}
GA_1(G)\geq \min\big\{GA_1(H_{\delta,\Delta}),GA_1(K_{\delta,\Delta})\big\}
.
\end{equation}
The following example shows that the answer is negative.

\begin{example}\label{ex:4-56} Let us suppose $\delta=4$ and $\Delta=56$. Consider a graph $G$ with $57$ vertices, two of them, $a_1,a_2$ with degree $56$ and the rest, $b_1,\dots,b_{55}$ with degree $4$. Let us assume the edges are as follows. There is an edge $a_ib_j$ for every $i,j$, an edge $a_1a_2$ and the vertices $b_1,\dots,b_{55}$ induce a cycle of length $55$. Note that these edges produce the claimed degree in each vertex.

Notice that $G$ has $166$ edges, one of them joins two vertices of degree $56$, $110$ of them join vertices with degree $56$ with vertices with degree $4$ and $55$ of them join vertices with degree $4$. Therefore,
$GA_1(G)=\frac{2\cdot 110 \sqrt{4\cdot 56}}{4+56} + 56=\frac{ 220 \sqrt{224}}{60} + 56\approx 110.8776.$

However,
$GA_1(H_{4,56})=\frac{112\sqrt{224}}{60}+84\approx 111.9377,$
and
$GA_1(K_{4,56})=\frac{448 \sqrt{224}}{60}\approx 111.7508. $

Also, by Proposition \ref{prop edges}, any graph with minimum degree $4$ and maximum degree $56$ has at least $140$ edges (while $G$ has $166$). By Theorem \ref{Th: lower bound} we have that any graph $G'$ with minimum degree $4$ and maximum degree $56$  satisfies that
$
GA_1(G')\geq  56 \cdot 5 \frac{ \sqrt{224}}{60}\approx 69.8443.
$
Notice that his lower bound is relatively far from the results obtained from $G$, $H_{4,56}$ and  $K_{4,56}$.
\end{example}

However, \eqref{eq:conj} holds if we have either $\d=1$ or $\d=\D$ (see Theorem \ref{Th: lower bound}).
Furthermore, \eqref{eq:conj} also holds if $\d$ and $\D$ are close enough, as
the following results show.

\begin{theorem}\label{Th:d-D-close}
Let $G$ be a graph with minimum degree $\delta>0$ and maximum degree $\Delta \ge 2$. If
\begin{equation}\label{Close h1}
\frac{2\sqrt{\delta\Delta}}{\delta+\Delta}\geq \frac{\Delta (\delta-1)}{\Delta (\delta-1)+2},
\end{equation}
then
\begin{equation}\label{Close 1}
GA_1(G) \geq \frac{2\Delta\sqrt{\delta\Delta}}{\delta+\Delta} + \frac{\Delta(\delta-1)}{2}.
\end{equation}

Furthermore, if $\Delta(\delta+1)$ is odd,
\begin{equation}\label{Close h2}
\frac{2\sqrt{\delta\Delta}}{\delta+\Delta}\geq \frac{\Delta (\delta-1)}{\Delta (\delta-1)+2} \ \mbox{ and } \ \frac{3\sqrt{\delta\Delta}}{\delta+\Delta}+\delta-\frac{1}{2}\geq \frac{2\sqrt{(\delta+1)\Delta}}{\delta+1+\Delta}+\frac{2\delta\sqrt{\delta(\delta+1)}}{2\delta+1},
\end{equation}
then
\begin{equation}\label{Close 2}
GA_1(G) \geq \frac{2(\Delta-1)\sqrt{\delta\Delta}}{\delta+\Delta} + \frac{2\sqrt{(\delta+1)\Delta}}{\delta+1+\Delta} +
\frac{2\delta \sqrt{\delta(\delta+1)}}{2\delta+1}+ \frac{(\Delta-2)(\delta-1)-1}{2}.
\end{equation}

If $\D$ and $\d$ verify \eqref{Close h1}, then the equality in \eqref{Close 1} is attained if and only if $\Delta(\delta+1)$ is even and $G \in \mathcal{G}_{\d,\D}$.
If $\D$ and $\d$ verify \eqref{Close h2} and $\Delta(\delta+1)$ is odd,
then the equality in \eqref{Close 2} is attained if and only if $G \in \mathcal{G}_{\d,\D}$.
\end{theorem}

\begin{proof}
Suppose $G\notin \mathcal{G}_{\d,\D}$. Then, by Proposition \ref{prop edges}, $G$ has at least $\frac{\Delta(\delta+1)}{2}+1$ edges. By (\ref{eq: bound}), this implies that
\begin{equation}
\label{equa1}
GA_1(G)\geq \Big(\frac{\Delta(\delta+1)}{2}+1\Big) \frac{2\sqrt{\delta \Delta}}{\delta+\Delta}.
\end{equation}
Let us denote $\varepsilon=\frac{2\sqrt{\delta \Delta}}{\delta+\Delta}$. Then, it suffices to check that
$$
\Big(\frac{\Delta(\delta+1)}{2}+1\Big) \varepsilon \geq \Delta \varepsilon +\frac{\Delta(\delta-1)}{2}.
$$
Thus, it is readily seen that
$$\Big(\frac{\Delta(\delta+1)}{2}+1\Big) \varepsilon \geq \Delta \varepsilon +\frac{\Delta(\delta-1)}{2} \quad \Leftrightarrow \quad
\varepsilon \Big(\frac{\Delta(\delta+1)}{2}+1-\Delta \Big)  \geq \frac{\Delta(\delta-1)}{2} \quad \Leftrightarrow $$
$$\varepsilon \Big(\frac{\Delta(\delta-1)+2}{2} \Big)  \geq \frac{\Delta(\delta-1)}{2} \quad \Leftrightarrow \quad
\varepsilon \geq \frac{\Delta(\delta-1)}{\Delta(\delta-1)+2}.$$

If $\Delta(\delta+1)$ is odd and $G\notin \mathcal{G}_{\d,\D}$, then by Proposition \ref{prop edges}, $G$ has at least
$\frac{\Delta(\delta+1)+1}{2}+1$ edges. By (\ref{eq: bound}), this implies that
\begin{equation}
\label{equa2}
GA_1(G)\geq \Big(\frac{\Delta(\delta+1)+1}{2}+1\Big) \frac{2\sqrt{\delta \Delta}}{\delta+\Delta}.
\end{equation}
Then, it suffices to check that
$$
\Big(\frac{\Delta(\delta+1)+1}{2}+1\Big) \varepsilon \ge (\Delta-1)\varepsilon + \frac{2\sqrt{(\delta+1)\Delta}}{\delta+1+\Delta} +
\frac{2\delta \sqrt{\delta(\delta+1)}}{2\delta+1}+ \frac{(\Delta-2)(\delta-1)-1}{2}.
$$
Since $\varepsilon \geq \frac{\Delta (\delta-1)}{\Delta (\delta-1)+2},$ the  argument from the even case implies that
$$\Big(\frac{\Delta(\delta+1)}{2}+1\Big) \varepsilon \geq \Delta \varepsilon +\frac{\Delta(\delta-1)}{2}=(\Delta-1) \varepsilon + \varepsilon +\frac{\Delta(\delta-1)}{2}$$
and it suffices to check that
$$\frac{3\varepsilon}{2}+\frac{\Delta(\delta-1)}{2}\geq \frac{2\sqrt{(\delta+1)\Delta}}{\delta+1+\Delta} +
\frac{2\delta \sqrt{\delta(\delta+1)}}{2\delta+1}+ \frac{(\Delta-2)(\delta-1)-1}{2}.$$
Since, by hypothesis
$$\frac{3}{2}\varepsilon+\delta-\frac{1}{2}\geq \frac{2\sqrt{(\delta+1)\Delta}}{\delta+1+\Delta}+\frac{2\delta\sqrt{\delta(\delta+1)}}{2\delta+1},$$
then the result follows immediately.

Proposition \ref{Prop minimal} gives that the equality in \eqref{Close 1} is attained if $\Delta(\delta+1)$ is even and $G \in \mathcal{G}_{\d,\D}$,
and that the equality in \eqref{Close 2} is attained if $\Delta(\delta+1)$ is odd and $G \in \mathcal{G}_{\d,\D}$.

Assume that $\D$ and $\d$ verify \eqref{Close h1}.
Proposition \ref{Prop minimal} gives that if the equality is attained in \eqref{Close 1} for some $G$, then $\Delta(\delta+1)$ is even.
Assume that the equality is attained in \eqref{Close 1} for some $G \notin \mathcal{G}_{\d,\D}$.
Thus, the equality is attained in \eqref{equa1}.
Remark \ref{remark bipartite} and \eqref{eq: bound} give that
$
|E(G)| = \frac{\Delta(\delta+1)}{2}+1
$
and $G$ is either regular or bipartite with the two sets being respectively the set of vertices with degree $\delta$ and degree $\Delta$.
If $G$ is regular, then $\d=\D$ and
$
\frac{\Delta(\Delta+1)}{2}+1 = |E(G)| = \frac{n\Delta}{2},
$
where $n=|V(G)|$.
So, $2=\D(n-1-\D)$ and, since $\Delta \ge 2$, we conclude $\D=2$ and $n-1-\D=1$; hence, $G$ is a $2$-regular graph with $n=4$ vertices, i.e., $G\cong C_4$.
Assume now that $G$ is bipartite with the two sets being respectively the set of vertices with degree $\delta$ and degree $\Delta$ (thus, $\D>\d$).
Then there exists a vertex $v_0\in V(G)$ of degree $\D$ with neighbors $v_1,\dots,v_\D\in V(G)$ of degree $\d$; since $G$ is a bipartite graph,
$$
\frac{\Delta(\delta+1)}{2}+1 = |E(G)| \ge \Delta\d,
\qquad
\Delta+2 \ge \Delta\d.
$$
If $\d \ge 2$, then $\Delta+2 \ge 2\Delta$, $2 \ge \Delta$ and we conclude $\Delta=\delta=2$, a contradiction.
If $\d = 1$, then $G$ is isomorphic to the union of $r$ graphs $K_{1,\D}$, since $G$ is bipartite, and we have
$$
\Delta +1 = \frac{\Delta(\delta+1)}{2}+1 = |E(G)| = r \Delta,
\qquad
1=(r-1) \Delta.
$$
Hence, $\Delta=1$, a contradiction.
We conclude that if the equality is attained in \eqref{Close 1} for some $G \notin \mathcal{G}_{\d,\D}$, then $(\D,\d)=(2,2)$ and $G\cong C_4$,
but we have
$
GA_1(C_4) =4 >3 = \frac{2\Delta\sqrt{\delta\Delta}}{\delta+\Delta} + \frac{\Delta(\delta-1)}{2},
$
a contradiction.
Therefore, if $G \notin \mathcal{G}_{\d,\D}$, then the inequality \eqref{Close 1} is strict.

Assume that $\D$ and $\d$ verify \eqref{Close h2}
and that the equality is attained in \eqref{Close 2} for some $G \notin \mathcal{G}_{\d,\D}$.
Thus, the equality is attained in \eqref{equa2}.
Remark \ref{remark bipartite} and \eqref{eq: bound} give that
$
|E(G)| = \frac{\Delta(\delta+1)+1}{2}+1
$
and $G$ is either regular or bipartite with the two sets being respectively the set of vertices with degree $\delta$ and degree $\Delta$.
If $G$ is regular, then $\d=\D$ and $\Delta(\Delta+1)$ is even, leading to contradiction with the number of edges.
%
Assume now that $G$ is bipartite with the two sets being respectively the set of vertices with degree $\delta$ and degree $\Delta$ (thus, $\D>\d$).
Then there exists a vertex $v_0\in V(G)$ of degree $\D$ with neighbors $v_1,\dots,v_\D\in V(G)$ of degree $\d$; since $G$ is a bipartite graph,
$$
\frac{\Delta(\delta+1)+1}{2}+1 = |E(G)| \ge \Delta\d,
\qquad
\Delta+3 \ge \Delta\d.
$$
Since $\Delta(\delta+1)$ is odd, we have $\d \ge 2$.
Thus, $\Delta+3 \ge 2\Delta$, $3 \ge \Delta > \d \ge 2$ and we conclude $\Delta=3$, $\delta=2$ and $|E(G)|=6$.
Hence, $G$ has two vertices with degree 3 and three vertices with degree 2 and it is isomorphic either to the cycle graph $C_5$ with an additional edge or to the complete bipartite graph $K_{2,3}$.
One can check that for these graphs the equality in \eqref{Close 2} is strict, a contradiction.
We conclude that the equality is not attained in \eqref{Close 2} if $G \notin \mathcal{G}_{\d,\D}$.
\end{proof}

%

\begin{remark}\label{remark: odd} Notice that in Theorem \ref{Th:d-D-close}, since $\frac{2\sqrt{(\delta+1)\Delta}}{\delta+1+\Delta}\leq 1$ for every $\delta<\Delta$, to assure the second condition in the case where
$\Delta(\delta+1)$ is odd, this is,
$$\frac{3\sqrt{\delta\Delta}}{\delta+\Delta}+\delta-\frac{1}{2}\geq \frac{2\sqrt{(\delta+1)\Delta}}{\delta+1+\Delta}+\frac{2\delta\sqrt{\delta(\delta+1)}}{2\delta+1},$$
it suffices to check that
$$\frac{3\sqrt{\delta\Delta}}{\delta+\Delta}+\delta-\frac{3}{2}\geq  \frac{2\delta\sqrt{\delta(\delta+1)}}{2\delta+1},$$
or equivalently,
$$\delta\Big(1- \frac{2\sqrt{\delta(\delta+1)}}{2\delta+1}\Big)\geq \frac{3}{2}\Big(1-\frac{2\sqrt{\delta\Delta}}{\delta+\Delta}\Big) .$$
\end{remark}

A nontrivial connected graph with maximum degree at most four is a \emph{molecular graph} representing hydrocarbons \cite{Trinajstic}.
Theorem \ref{Th:d-D-close} allows to obtain sharp inequalities for molecular graphs.

\begin{corollary}\label{C:d-hd 4}
Let $G$ be a molecular graph with minimum degree $\delta>0$ and maximum degree $\Delta$.
If $(\delta,\Delta) \neq (2,3)$, then
$$
GA_1(G)
\geq \frac{2\Delta\sqrt{\delta\Delta}}{\d+\Delta} + \frac{\Delta(\delta-1)}{2},
$$
with equality if and only if $G \in \mathcal{G}_{\d,\D}$.
If $(\delta,\Delta) = (2,3)$, then
$$
GA_1(G)
\geq \frac{2(\Delta-1)\sqrt{\delta\Delta}}{\delta+\Delta} + \frac{2\sqrt{(\delta+1)\Delta}}{\delta+1+\Delta} +
\frac{2\delta \sqrt{\delta(\delta+1)}}{2\delta+1}+ \frac{(\Delta-2)(\delta-1)-1}{2}
=\frac{8\sqrt{6}}{5}+ 1,
$$
with equality if and only if $G \in \mathcal{G}_{2,3}$.
\end{corollary}

\begin{proof}
Since $G$ is a molecular graph, $\Delta(\delta+1)$ is even if and only if $(\delta,\Delta) \neq (2,3)$.
One can check that, in this case, $(\delta,\Delta)$ satisfies \eqref{Close h1} in Theorem \ref{Th:d-D-close}.
Thus, Theorem \ref{Th:d-D-close} gives the first part of the corollary.
It is easy to check that $(\delta,\Delta) = (2,3)$ satisfies \eqref{Close h2} in Theorem \ref{Th:d-D-close}.
Hence, Theorem \ref{Th:d-D-close} gives the second part of the corollary.
\end{proof}

\begin{corollary}\label{C:d-hd}
Let $G$ be a graph with  minimum degree $\delta>0$ and maximum degree $\Delta=\delta+h \ge 2$. If
$(16-h^2)\Delta^3+(2h^3+2h^2-32h-16)\Delta^2+(-h^4-2h^3+15h^2+16h+16)\Delta-16h\geq 0,$
then
$$
\begin{aligned}
GA_1(G) &
\geq \frac{2\Delta\sqrt{\Delta(\Delta-h)}}{2\Delta-h} + \frac{\Delta(\Delta-h-1)}{2}.
\end{aligned}
$$
\end{corollary}

\begin{proof} By Theorem \ref{Th:d-D-close} and $\delta=\D-h$, it suffices to check that
$$\frac{2\sqrt{\Delta(\Delta-h)}}{2\Delta-h}\geq \frac{\Delta (\Delta-(h+1))}{\Delta (\Delta-(h+1))+2}.$$

Thus,
$$
\begin{aligned}
& \frac{2\sqrt{\Delta(\Delta-h)}}{2\Delta-h}\geq \frac{\Delta (\Delta-(h+1))}{\Delta (\Delta-(h+1))+2}  \Leftrightarrow
2\sqrt{\Delta(\Delta-h)}\geq \frac{\Delta (\Delta-(h+1))(2\Delta-h)}{\Delta (\Delta-(h+1))+2}  \Leftrightarrow
\\
& 4(\Delta -h)\geq \frac{\Delta(\Delta-(h+1))^2(2\Delta-h)^2}{\Delta^2 (\Delta-(h+1))^2+4\Delta(\Delta-(h+1))+4}  \Leftrightarrow
\\
& 4(\Delta -h)\Delta^2 (\Delta-(h+1))^2+ 16\Delta(\Delta -h)(\Delta-(h+1))+16(\Delta -h) 
\\
& \qquad \qquad \geq \Delta(\Delta-(h+1))^2(2\Delta-h)^2
\Leftrightarrow
\\
& \Delta (\Delta-(h+1))^2[4\Delta(\Delta -h)-(2\Delta-h)^2]+16(\Delta -h)[\Delta(\Delta-(h+1))+1] \geq 0
 \Leftrightarrow
\\
& -h^2\Delta (\Delta-(h+1))^2+16(\Delta -h)[\Delta^2-h\Delta-\Delta +1] \geq 0  \Leftrightarrow
\\
& (16-h^2)\Delta^3+(2h^3+2h^2-32h-16)\Delta^2+(-h^4-2h^3+15h^2+16h+16)\Delta-16h \geq 0.
\end{aligned}
$$
\end{proof}

Let us denote $$P(h,\Delta)=(16-h^2)\Delta^3+(2h^3+2h^2-32h-16)\Delta^2+(-h^4-2h^3+15h^2+16h+16)\Delta-16h.$$
Therefore, we obtain the following polynomials with the following real solutions (rounded off to one decimal):

If $h=0$, $P(0,\Delta)=16\Delta^3 -16 \Delta^2 + 16\Delta$, real root: $0$.

If $h=1$, $P(1,\Delta)=15\Delta^3 -44 \Delta^2 + 44\Delta -16$, real root: $1.3$.

If $h=2$, $P(2,\Delta)=12\Delta^3 -56 \Delta^2 + 76\Delta -32 $, real root: $2.7$.

If $h=3$, $P(3,\Delta)=7\Delta^3 - 40\Delta^2 + 64\Delta - 48$, real root: $3.8$.

If $h=4$, $P(4,\Delta)=16 \Delta^2 - 64 \Delta - 64$, real roots: $-0.8$ and $4.8$.

If $h=5$, $P(5,\Delta)=-9\Delta^3 + 124\Delta^2 -404 \Delta -80 $, real roots: $-0.2$, $5.9$ and $8.1$.

If $h=6$, $P(6,\Delta)=-20 \Delta^3 + 296\Delta^2 -1076 \Delta - 96$, real roots: $-0.1$, $6.9$ and $8$.

If $h=7$, $P(7,\Delta)=-33\Delta^3 + 544\Delta^2 -2224 \Delta - 112$, real roots: $-0.1$, $7.9$ and $8.6$.

This, together with Theorem \ref{Th: lower bound} and Corollary \ref{C:d-hd}, yields the following:

\begin{corollary}\label{C:d-hd 2}
Let $G$ be a graph with  minimum degree $\delta>0$ and maximum degree $\Delta=\delta+h \ge 2$.
If we have

$(1)$ $h=0$ or $h=1$, for every $\Delta\geq 2$,

$(2)$ $h=2$, for every $\Delta\geq 3$,

$(3)$ $h=3$, for every $\Delta\geq 4$,

$(4)$ $h=4$, for every $\Delta\geq 5$,

$(5)$ $h=5$, for every $\Delta\in \{6,7,8\}$,

$(6)$ $h=6$, for every $\Delta\in \{7,8\}$,

$(7)$ $h\geq 7$ and $\Delta= h+1$,

\noindent then
$$
\begin{aligned}
GA_1(G) &
\geq \frac{2\Delta\sqrt{\Delta(\Delta-h)}}{2\Delta-h} + \frac{\Delta(\Delta-h-1)}{2}.
\end{aligned}
$$
\end{corollary}

\begin{corollary}\label{C:d-1d}
Let $G$ be a graph with maximum degree $\Delta \ge 2$ and minimum degree $\delta=\D-1$. Then
$$
\begin{aligned}
GA_1(G) &
\geq \frac{2\Delta\sqrt{\Delta(\Delta-1)}}{2\Delta-1} + \frac{\Delta(\Delta-2)}{2}, \quad \text{if $\Delta$ is even,}
\\
GA_1(G) &
\geq \frac{4(\Delta-1)\sqrt{\Delta(\Delta-1)}}{2\Delta-1} + \frac{(\Delta-2)^2-1}{2}+1, \quad \text{if $\Delta$ is odd,}
\end{aligned}
$$
with equalities if and only if $G \in \mathcal{G}_{\D-1,\D}$.
\end{corollary}

\begin{proof} 
%
As we saw above, $P(1,\D)=15\Delta^3-44\Delta^2+44\Delta-16\geq 0$
for every $\Delta \geq 2$. Furthermore, since this inequality is strict for every $\Delta\geq 2$, the bound is only
attained if $G\in \mathcal{G}_{\Delta-1,\Delta}$.

If $\D$ is odd (and therefore, $\D\geq 3$), then by Theorem \ref{Th:d-D-close}, and Remark \ref{remark: odd}, the second result follows trivially from the fact that $\delta\geq 2>\frac{3}{2}$.
Also, since the inequality is strict, the bound is only attained if $G\in \mathcal{G}_{\Delta-1,\Delta}$
\end{proof}

Corollary \ref{C:d-hd 2} has also the following consequence.

\begin{corollary}\label{C:d-hd 3}
Let $G$ be a graph with  minimum degree $\delta>0$ and maximum degree $2 \le \Delta \le 8$.
Then
$$
\begin{aligned}
GA_1(G) &
\geq \frac{2\Delta\sqrt{\delta\Delta}}{\d+\Delta} + \frac{\Delta(\delta-1)}{2}.
\end{aligned}
$$
\end{corollary}

\begin{lemma}\label{max-deg} If $28\leq a \leq \Delta$, then
\begin{equation}\label{eq:max-deg}
\frac{2\sqrt{a\Delta}}{\Delta+a}> 2 \frac{2\sqrt{2\Delta}}{\Delta+2}.
\end{equation}
\end{lemma}

\begin{proof}
$$
\begin{aligned}
\frac{2\sqrt{a\Delta}}{\Delta+a}> 2 \frac{2\sqrt{2\Delta}}{\Delta+2} \quad \Leftrightarrow \quad
\frac{\sqrt{a}}{\Delta+a}>  \frac{2\sqrt{2}}{\Delta+2} \quad \Leftrightarrow
\\
a(\Delta+2)^2> 8(\Delta+a)^2 \quad \Leftrightarrow \quad
(a-8)\Delta^2-12a\Delta-8a^2+4a > 0.
\end{aligned}
$$
Since $28\leq a \leq \Delta$, then $(a-8)\Delta^2-12a\Delta-8a^2+4a \geq (a-28)\Delta^2+4a>0$ and (\ref{eq:max-deg}) holds.
\end{proof}

\begin{lemma}\label{min-deg} If $2\leq b\leq 27$ and $\Delta\geq 30$, then
\begin{equation}\label{eq:min-deg}
(b-1)\frac{2\sqrt{2b}}{b+2}> b \frac{2\sqrt{2\Delta}}{\Delta+2}.
\end{equation}
\end{lemma}

\begin{proof}
We have
$$
(b-1)\frac{2\sqrt{2b}}{b+2}> b \frac{2\sqrt{2\Delta}}{\Delta+2}
\quad
\Leftrightarrow
\quad
\frac{b-1}{(b+2)\sqrt{b}}>  \frac{\sqrt{\Delta}}{\Delta+2} .
$$
Let us define
$$
A(b)=\frac{b-1}{(b+2)\sqrt{b}},
\qquad
B(\D)= \frac{\sqrt{\Delta}}{\Delta+2} .
$$
One can check that
$\min_{2\leq b\leq 27}A(b)=A(27)$ and
$\max_{\D \ge 30} B(\D)= B(30)$.
Since
$
A(27)
= \frac{26}{29\sqrt{27}}
> \frac{\sqrt{30}}{32} = B(30) ,
$
we conclude that (\ref{eq:min-deg}) holds.
\end{proof}

\begin{theorem} \label{l:230} Let $G$ be a graph with  minimum degree $2$ and maximum degree $\Delta\geq 28$. Then,
\[GA_1(G)\geq 2\Delta \frac{2\sqrt{2\Delta}}{\Delta+2}, \]
and the equality is attained if and only if $G=K_{2,\Delta}$.
\end{theorem}

\begin{proof} Let $x_0$ be a vertex such that $d_{x_0}=\Delta$. Let $C_1,\dots,C_k$ be the connected components of $G\backslash \{x_0\}$ and
$R_i=V(C_i)\cap N(x_0)$ for every $1\leq i \leq k$ where $N(x_0)$ denotes the set of vertices in $G$ adjacent to $x_0.$
Let $r_i:=|R_i|$ and notice that $\sum_{i=1}^kr_i=\Delta$.
Denote by $\G_i$ the subgraph of $G$ induced by $V(C_i)\cup \{x_0\}$ (thus, $\cup_{i=1}^k \G_i=G$).

If $C_i$ has at least $r_i$ edges, then $\sum_{uv\in E(C_i)}\frac{2\sqrt{d_ud_v}}{d_u+d_v}\geq r_i \frac{2\sqrt{2\Delta}}{\Delta+2}$.
Since $|E(\G_i)| = |E(C_i)| +r_i$, we have that
\[\sum_{uv\in E(\G_i)}\frac{2\sqrt{d_ud_v}}{d_u+d_v} \geq 2r_i \frac{2\sqrt{2\Delta}}{\Delta+2}.\]

If $C_i$ has less than $r_i$ edges then, since $R_i$ has $r_i$ vertices and $R_i\subseteq C_i$ with $C_i$ connected, it follows that $R_i=V(C_i)$ and $C_i$ is a tree with exactly $r_i-1$ edges.

Suppose there is a vertex $v\in R_i$, such that $d_v\geq 28$.
Recall that $v$ is adjacent to $x_0$ with $d_{x_0}=\Delta$.
Thus, by Lemma \ref{max-deg}, $\frac{2\sqrt{d_v\Delta}}{\Delta+d_v}> 2 \frac{2\sqrt{2\Delta}}{\Delta+2}$.
Since apart from the edge $x_0v$, there are $r_i-1$ edges joining $x_0$ and the vertices in $R_i\backslash \{v\}$ and $r_i-1$ edges in $C_i$, it follows that
\[\sum_{uv\in E(\G_i)}\frac{2\sqrt{d_ud_v}}{d_u+d_v}> 2r_i \frac{2\sqrt{2\Delta}}{\Delta+2}.\]

Otherwise, if $v_i$ is the vertex in $R_i$ with maximum degree and $d_{v_i} \le 27$, then
$
\frac{2\sqrt{d_ud_{v_i}}}{d_u+d_{v_i}}
\ge
\frac{2\sqrt{2d_{v_i}}}{d_{v_i}+2}
$
for every vertex $u\in N(v_i)\setminus \{x_0\}$. Therefore, by  Lemma \ref{min-deg}, if $A(d_{v_i})>B(\D)$
$$
\sum_{u\in N(v_i)\setminus \{x_0\}} \frac{2\sqrt{d_ud_{v_i}}}{d_u+d_{v_i}}
\ge
(d_{v_i}-1) \frac{2\sqrt{2d_{v_i}}}{d_{v_i}+2}
> d_{v_i} \frac{2\sqrt{2\D}}{\D+2}.
$$
Since apart from these edges there are $r_i-d_{v_i}$ edges in $C_i$ and
$r_i$ edges joining $x_0$ and the vertices in $R_i$, it follows that
\[\sum_{uv\in E(\G_i)}\frac{2\sqrt{d_ud_v}}{d_u+d_v}> 2r_i \frac{2\sqrt{2\Delta}}{\Delta+2}.\]

Notice that $A(d_{v_i})>B(\D)$ for every $\D\geq 30$, if $\D=29$ and $d_{v_i}\leq 26$ and if $\D=28$ and $d_{v_i}\leq 25$.
Therefore, if every component $C_i$ either satisfies one of these cases, or has $r_i$ edges, or has a vertex with degree at least 28, then
$$GA_1(G)=\sum_{i=1}^k \sum_{uv\in E(\G_i)} \frac{2\sqrt{d_ud_v}}{d_u+d_v} >
\sum_{i=1}^k 2r_i \frac{2\sqrt{2\Delta}}{\Delta+2}
=2\Delta \frac{2\sqrt{2\Delta}}{\Delta+2} =GA_1(K_{2,\Delta}). $$

Therefore, to finish the proof it suffices to check the following cases:

Case 1. Suppose $\D=29$ and there is a vertex $v_i$ adjacent to $x_0$ such that $d_{v_i}= 27$, $v_i\in C_1$ and $C_1$ has $r_1-1$ edges.
Then, there are exactly two vertices adjacent to $x_0$ which are not adjacent to $v_i$.
Let us assume, relabeling if necessary, that these are $x_{28}$ and $x_{29}$, and $v_i=x_1$.
Therefore, $G$ has one edge joining $x_0$ to $x_1$ where $d_{x_0}=29$ and $d_{x_1}=27$, $28$ edges joining $x_0$ to a vertex $x_i$ with $d_{x_i}\geq 2$ and $26$ edges joining $x_1$ to $x_j$ for $2\leq j \leq 27$ with $d_{x_j}\leq 27$.
If $G$ has $58$ edges, then trivially $GA_1(G)\ge GA_1(K_{2,29})$. If
$|E(G)|\leq 57$, then there are at most two edges left.
Since $d_{x_{28}}\geq 2$ and $d_{x_{29}}\geq 2$ either there is an edge $x_{28}x_{29}$ or there are two edges joining $x_{28}$ and $x_{29}$ to the same or two different vertices in $\{x_2,\dots,x_{27}\}$ or there is an edge $x_{28}x_{29}$ and some extra edge joining two vertices in
$\{x_2,\dots,x_{29}\}$. Thus, either there is an edge joining two vertices with degree $2$ or two edges joining a vertex with degree $2$ or $3$ to a vertex with degree $3$ or $4$. Since
$2\frac{2\sqrt{2\cdot 3}}{2+3}>2\frac{2\sqrt{2\cdot 4}}{2+4}>1$,
it follows that
$$GA_1(G)\geq \frac{2\sqrt{29\cdot 27}}{56} + 28 \frac{2\sqrt{29\cdot 2}}{31} + 26\frac{2\sqrt{27\cdot 2}}{29}+1$$
and it suffices to check that $1.999\approx \frac{2\sqrt{29\cdot 27}}{56}+1> 4\frac{2\sqrt{29\cdot 2}}{31}\approx 1.965$.

Case 2. Suppose $\D =28$ and there is a vertex $v_i$ adjacent to $x_0$ such that $d_{v_i}= 26$, $v_i\in C_1$ and $C_1$ has $r_1-1$ edges.
Then, there are exactly two vertices adjacent to $x_0$ which are not adjacent to $v_i$.
Let us assume, relabeling if necessary, that these are $x_{27}$ and $x_{28}$, and $v_i=x_1$.
Therefore, $G$ has one edge joining $x_0$ to $x_1$ where $d_{x_0}=28$ and $d_{x_1}=26$,
$27$ edges joining $x_0$ to a vertex $x_i$ with $d_{x_i}\geq 2$ and
$25$ edges joining $x_1$ to $x_j$ for $2\leq j \leq 26$ with $d_{x_j}\leq 26$.
If $G$ has $56$ edges, then trivially $GA_1(G)\ge GA_1(K_{2,28})$.
If $|E(G)|\leq 55$, then there are at most two edges left.
Since $d_{x_{27}}\geq 2$ and $d_{x_{28}}\geq 2$ either there is an edge $x_{27}x_{28}$ or there are two edges joining $x_{27}$ and $x_{28}$ to the same or two different vertices in $\{x_2,\dots,x_{26}\}$ or there is an edge $x_{27}x_{28}$ and some extra edge joining two vertices in $\{x_2,\dots,x_{28}\}$. Thus, either there is an edge joining two vertices with degree $2$ or two edges joining a vertex with degree $2$ or $3$ to a vertex with degree $3$ or $4$. Since
$2\frac{2\sqrt{2\cdot 3}}{2+3}>2\frac{2\sqrt{2\cdot 4}}{2+4}>1$,
it follows that
$$GA_1(G)> \frac{2\sqrt{28\cdot 26}}{54} + 27 \frac{2\sqrt{28\cdot 2}}{30} + 25\frac{2\sqrt{26\cdot 2}}{28}+1$$
and it suffices to check that $1.999\approx \frac{2\sqrt{28\cdot 26}}{54}+1> 4\frac{2\sqrt{28\cdot 2}}{30}\approx 1.996$.

Case 3. Suppose $\D =28$ and there is a vertex $v_i$ adjacent to $x_0$ such that $d_{v_i}= 27$, $v_i\in C_1$ and $C_1$ has $r_1-1$ edges.
Then, there is exactly one vertex adjacent to $x_0$ which is not adjacent to $v_i$.
Let us assume, that it is $x_{28}$, and $v_i=x_1$.
Therefore, $G$ has one edge joining $x_0$ to $x_1$ where $d_{x_0}=28$ and $d_{x_1}=27$,
$27$ edges joining $x_0$ to a vertex $x_i$ with $d_{x_i}\geq 2$ and
$26$ edges joining $x_1$ to $x_j$ for $2\leq j \leq 27$ with $d_{x_j}\leq 27$.
If $G$ has $56$ edges, then trivially $GA_1(G)\ge GA_1(K_{2,28})$.
If $|E(G)|\leq 55$, then there is at most 1 edge left.
Since $d_{x_{28}}\geq 2$ there is an edge joining $x_{28}$  to a vertex in $\{x_2,\dots,x_{27}\}$.
Thus, there is an edge joining a vertex with degree $2$ to a vertex with degree $3$.
Hence, it follows that
$$GA_1(G)\geq \frac{2\sqrt{28\cdot 27}}{55} + 27 \frac{2\sqrt{28\cdot 2}}{30} + 26\frac{2\sqrt{27\cdot 2}}{29}+\frac{2\sqrt{2\cdot 3}}{5}$$
and it suffices to check that
$1.98 \approx \frac{2\sqrt{28\cdot 27}}{55}+\frac{2\sqrt{2\cdot 3}}{5}> 3\frac{2\sqrt{28\cdot 2}}{30}\approx 1.50$.

Since the inequalities in Lemmas \ref{max-deg} and \ref{min-deg} are strict, it follows from the proof that if $GA_1(G)=2\Delta \frac{2\sqrt{2\Delta}}{\Delta+2}$,
then $C_i$ has exactly $r_i$ edges for every $1\leq i \leq k$.
Therefore, $G$ has $2\Delta$ edges and Proposition \ref{prop vertex-edges} gives $G=K_{2,\Delta}$.
\end{proof}

Given any odd integer $\Delta\geq 3$, let us define $H_{\D}$ as the  graph  with minimum degree $2$, maximum degree $\Delta$, $|V(H_\D)|=\D+1$, and such that there are $2$ vertices, $x_0, x_1$ with degree $\D$ which are adjacent and $\D-1$ vertices with degree $2$: $x_2,\dots,x_{\D}$.
Note that
\begin{equation}\label{eq: H}
GA_1(H_\D) =2(\Delta-1)\frac{2\sqrt{2 \D}}{2+ \D}+1.
\end{equation}

The next result shows that the conclusion of Theorem \ref{l:230} does not hold for $\D< 28$.

\begin{proposition} For any integer $2\leq \Delta\leq 27$, if $G\in \mathcal{G}_{2,\Delta}$, then
$$
\begin{aligned}
GA_1(G) &
< GA_1 (K_{2,\D}), \quad \text{if $\Delta$ is even,}
\\
GA_1(H_\D) &
< GA_1 (K_{2,\D}), \quad \text{if $\Delta$ is odd.}
\end{aligned}
$$
\end{proposition}

\begin{proof} If $\D$ is even, then $\D(\d+1)$ is even and, by Proposition \ref{Prop: comp}, $GA_1(G) < GA_1 (K_{2,\D})$.

If $\D$ is odd, then by (\ref{eq: H}),
$GA_1(H_\D)=2(\D-1)\frac{2\sqrt{2\D}}{2+\D}+1$
and, since for every $2\leq \D\leq 27$ we have
$1<2\frac{2\sqrt{2\D}}{2+\D},$
we conclude
$GA_1(H_\D)<GA_1(K_{2,\D}).$
\end{proof}


\begin{thebibliography}{99}

\bibitem{Das10b} K. C. Das, On geometric-arithmetic index of graphs, \emph{MATCH Commun. Math. Comput. Chem.} {\bf 64} (2010) 619--630.

\bibitem{DGF} K. C. Das, I. Gutman, B. Furtula, Survey on Geometric-Arithmetic Indices of Graphs,
\emph{MATCH Commun. Math. Comput. Chem.} {\bf 65} (2011) 595--644.

\bibitem{DGF2} K. C. Das, I. Gutman, B. Furtula, On first geometric-arithmetic index of graphs,
\emph{Discrete Appl. Math.} {\bf 159} (2011) 2030--2037.

\bibitem{GF} I. Gutman, B. Furtula (Eds.), Recent Results in the Theory of Randi\'c Index, Univ.
Kragujevac, Kragujevac, 2008.

\bibitem{LG} X. Li, I. Gutman, Mathematical Aspects of Randi\'c Type Molecular Structure Descriptors,
Univ. Kragujevac, Kragujevac, 2006.

\bibitem{LS} X. Li, Y. Shi, A survey on the Randi\'c index, \emph{MATCH Commun. Math. Comput. Chem.} {\bf 59} (2008) 127--156.

\bibitem{MH} M. Mogharrab, G. H. Fath-Tabar, Some bounds on $GA_1$ index of graphs,
\emph{MATCH Commun. Math. Comput. Chem.} {\bf 65} (2010) 33--38.

\bibitem{R} M. Randi\'c, On characterization of molecular branching, \emph{J. Am. Chem. Soc.} {\bf 97} (1975) 6609--6615.

\bibitem{RS2} J. M. Rodr{\'\i}guez, J. M. Sigarreta, On the Geometric-Arithmetic Index,
\emph{MATCH Commun. Math. Comput. Chem.} {\bf 74} (2015) 103--120.

\bibitem{RS3} J. M. Rodr{\'\i}guez, J. M. Sigarreta, Spectral properties of geometric-arithmetic index,
\emph{Appl. Math. Comput.} {\bf 277} (2016) 142--153.

\bibitem{RS} J. A. Rodr{\'\i}guez-Vel\'azquez, J. M. Sigarreta, On the Randi\'c index and condicional parameters of a graph,
\emph{MATCH Commun. Math. Comput. Chem.} {\bf 54} (2005) 403--416.

\bibitem{RS0} J. A. Rodr{\'\i}guez-Vel\'azquez, J. Tom\'as-Andreu, On the Randi\'c index of polymeric networks modelled by generalized Sierpinski graphs,
\emph{MATCH Commun. Math. Comput. Chem.} {\bf 74} (2015) 145--160.

\bibitem{S} J. M. Sigarreta, Bounds for the geometric-arithmetic index of a graph,
\emph{Miskolc Math. Notes} {\bf 16} (2015) 1199--1212.

\bibitem{TRC} TRC Thermodynamic Tables. Hydrocarbons; Thermodynamic Research Center,
The Texas A $\&$ M University System: College Station, TX, 1987.

\bibitem{Trinajstic}
N. Trinajsti\'c, Chemical Graph Theory, CRC Press, Boca Raton, Fl, USA, 2nd edition, 1992.

\bibitem{NGJ} M. V\"oge, A. J. Guttmann, I. Jensen, On the number of benzenoid hydrocarbons, \emph{J. Chem. Inf. Comput. Sci.} {\bf 42} (2002) 456--466.

\bibitem{VF} D. Vuki\v{c}evi\'c, B. Furtula, Topological index based on the ratios of geometrical and
arithmetical means of end-vertex degrees of edges, \emph{J. Math. Chem.} {\bf 46} (2009) 1369--1376.

\bibitem{Wi} H. Wiener, Structural determination of paraffin boiling points, {\it J. Am. Chem. Soc.} {\bf 69} (1947) 17--20.

\end{thebibliography}
\end{document}